\UseRawInputEncoding
\documentclass{article}

\usepackage{stmaryrd}
\usepackage{lmodern}
\usepackage{diagbox}
\usepackage{mathtools}
\usepackage{amssymb,amsmath,amsfonts}

\usepackage{mathrsfs}

\usepackage{color}
\usepackage{epsfig,graphics,graphicx}

\usepackage{lscape}
\usepackage{latexsym}

\def\proof{\noindent {\bf Proof }}
\def\alp{{\alpha}}

\def\gam{{\gamma}}
\def\gama{{\gamma}}
\def\qed{~~\vrule height8pt width4pt depth0pt}
\def\ex{{\mathbb E}}

\def\Xb{{\bar X}}
\def\Hb{{\bar H}}

\def\qb{{n}}

\def\Scal{{\cal S}}

\def\Mcal{{\cal M}}
\def\Ecal{{\cal E}}

\def\RS{{\cal RS}}

\newtheorem{prop}{Proposition}

\newtheorem{lemma}{Lemma}
\newtheorem{thm}{Theorem}
\newtheorem{cor}{Corollary}

\newcommand{\delO}{{\delta_0}}
\newcommand{\eps}{{\varepsilon}}

\newcommand{\fq}{{\mathbb F}_{q}}

\newcommand{\Ex}{{\mathbb E}}
\newcommand{\prob}{{\mathbb P}}
\newcommand{\Gh}{{\hat G}}

\begin{document}

\title {Asymptotic distributions of the number of zeros of random polynomials in Hayes equivalence class
over a finite field}
 \author{
Zhicheng Gao\\
School of Mathematics and Statistics\\
Carleton University\\
Ottawa, Ontario\\
Canada K1S5B6\\
Email:~zgao@math.carleton.ca
 }
\maketitle

\begin{abstract}
Hayes equivalence is defined on monic polynomials over a finite field $\fq$ in terms of the prescribed leading coefficients and the residue classes modulo a given monic polynomial $Q$. We study the distribution of the number of zeros in a random polynomial over finite fields in a given Hayes equivalence class. It is well known that the number of distinct zeros of a random polynomial over $\fq$ is asymptotically Poisson with mean 1. We show that this is also true for random polynomials in any given Hayes equivalence class. Asymptotic formulas are also given for the number of such polynomials when  the degree of such polynomials is proportional to $q$ and the degree of $Q$ and the number of prescribed leading coefficients are bounded by $\sqrt{q}$. When $Q=1$, the problem is equivalent to the study of the distance distribution in Reed-Solomon codes. Our asymptotic formulas extend some earlier results and imply that all words for a large family of Reed-Solomon codes are ordinary, which further supports the well-known {\em Deep-Hole} Conjecture.
\end{abstract}

\section{ Introduction}

Let $\fq$ be the finite field with $q$ elements of characteristic $p$. In this paper, we study the distribution of the number of zeros in a random polynomial over $\fq$ in a given Hayes equivalence class (the precise definition of Hayes equivalence will be given below). There are two motivations for the current study. First, there has been considerable interest in the study of distributions of parameters in general combinatorial structures.  See, e.g., \cite{ABT97,FS09} for general combinatorial structures, and \cite{GP06,KK90,KK93,KnoWal95,KW21,ZWW17} for polynomials over finite fields with respect to factorization patterns. Second, there is a close connection between the distance distribution over Reed-Solomon codes and the distribution of the number of distinct zeros in a random polynomial over $\fq$ in a given equivalence class defined by leading coefficients; see \cite{Gao21,GaoLi22,LiWan20,ZWW17} for related discussions.

Many interesting parameters of random combinatorial structures are known to be asymptotically Poisson \cite{ABT97,FS09}. Two well-known examples are the number of small cycles in a random permutation and the number of zeros in a random polynomial over $\fq$. There are basically two approaches to the study of such problems. One is based on analytic combinatorics \cite{FS09} using  generating functions, and the other is probabilistic using the famous Chen-Stein method \cite{AGG89,BHJ92,CDM05}. The latter approach requires that the random variable under consideration can be expressed as a sum of nearly independent indicator variables.

%The distribution of the number of components in general combinatorial structures with restricted patterns are studied in \cite{DGP09}.

In this paper, we apply the generating function approach to study the number of distinct zeros in a random polynomial over $\fq$ in a given Hayes equivalence class. This includes polynomials with prescribed leading and/or ending coefficients.

Before stating our main results, we  introduce some notations, which will be  used throughout the paper.
\begin{itemize}
\item $\Mcal$ denotes the set of monic polynomials over $\fq$, and $\Mcal_d$ denotes the set of polynomials in $\Mcal$ of degree $d$.
\item $\deg(f)$ denotes the degree of the polynomial $f$,  $[x^j]f$ denotes the coefficient of $x^j$ in $f$, and  ${\hat f}=x^{\deg(f)} f(1/x)$.
\end{itemize}

Fix a non-negative integer $\ell$ and a polynomial $Q\in \Mcal$.  Two polynomials $f,g \in \Mcal$ are said to be {\em Hayes equivalent} with respect to $\ell$ and $Q$ if
$\gcd(f,Q)=\gcd(g,Q)=1$ and
\begin{align}
{\hat f}(x)&\equiv {\hat g}(x) \pmod{x^{\ell+1}}, \label{eq:EL}\\
f(x)&\equiv g(x) \pmod{Q}. \label{eq:EQ}
\end{align}

The following two special cases are particularly interesting.

\begin{itemize}
\item[(a)] $Q=1$. In this case, condition \eqref{eq:EQ} is null, and Hayes equivalence is defined by the $\ell$ {\em  leading coefficients}  $[x^{\deg(f)-j}]f,1\le j\le \ell$.
\item[(b)] $Q=x^t$ for some $t>0$. In this case, Hayes equivalence is defined by the $\ell$ leading coefficients and $t$ {\em ending coefficients}  $[x^j]f, 0\le j< t$.
\end{itemize}

Let $\Ecal^{\ell,Q}$ denote the set of all Hayes equivalence classes with respect to $\ell,Q$, and let $\langle f\rangle$ denote the equivalence class represented by a polynomial $f\in \Mcal$. It is  known \cite{EH91,Hay65,Hsu96} that $\Ecal^{\ell,Q}$ is a group under the operation $\langle f\rangle \langle g\rangle=\langle fg\rangle$.

Given a set $D\subseteq \fq$ and $\eps\in \Ecal^{\ell,Q}$, let $\Mcal_{k+t+\ell}(\eps)$ denote the set of polynomials in $\Mcal_{k+t+\ell}$ which are equivalent to $\eps$, and  $Y_k(\eps)$ be the
number of zeros in $D$ of a random polynomial $f\in \Mcal_{k+t+\ell}(\eps)$ (under uniform distribution).

Since $\gcd(f,Q)=1$, {\em we will assume, without loss of generality, that $D$ does not contain any zero of $Q$.} In the rest of the paper, we will also set $n:=|D|$.

The paper \cite{Gao21} focuses on obtaining exact expression  of the distribution of $Y_k(\eps)$  with $Q=1$ and $\ell\le 2$. The paper \cite{GaoLi22} focuses on the problem whether $\prob(Y_k(\eps)=k+1)>0$ when $Q=1$ and $k$ is large, which is motivated by the {\em Deep-Hole Conjecture} about Reed-Solomon codes. In this paper we study the asymptotic distribution of $Y_k(\eps)$  for general $\ell$ and $Q$. 
 Our main results are summarized below.
\begin{thm}\label{cor:Binomial}
Let $Q\in \Mcal_t$ and  $\eps\in \Ecal^{\ell,Q}$.
\begin{itemize}
\item[(a)]  As $k-r\to \infty$ we have
\begin{align*}
\prob(Y_k(\eps)=r)&\sim {n\choose r}\left(\frac{1}{q}\right)^{r}\left(1-\frac{1}{q}\right)^{n-r}.
\end{align*}
That is, $Y_k(\eps)$ is asymptotically binomial.
\item[(b)]  
As $n,k-r\to \infty$ and for  $r=o(\sqrt{n})$, we have
\begin{align*}
\prob(Y_k(\eps)=r)&\sim e^{-n/q}\frac{1}{r!}\left(\frac{\qb}{q}\right)^r.
\end{align*}
That is, $Y_k(\eps)$ is asymptotically Poisson with mean $\qb/q$. 
\end{itemize}
\end{thm}

Recall that $p$ is the characteristic of $\fq$. Some applications (see e.g. \cite{GaoLi22}) require asymptotic formulas  of $\prob(Y_k(\eps)=r)$ when $r$ is close to $k$.  Let $c:=k/q$, $\gama:=(t+\ell-1)\sqrt{q}/n$, and $\delO$ denote any small positive constant (independent of $k$). The next  theorem provides a simple asymptotic formula for $\prob(Y_k(\eps)=r)$ which is uniform over all values of $r$  under  either of the following two conditions.
\begin{description}
\item[Condition A:] $\quad \displaystyle \frac{p-1}{p}c\ln \frac{1}{c}+(1-c)\ln \frac{1}{1-c}-\frac{1+c}{p}\ln(1+c)\ge \gama\ln(2p)+ \delO$.
\item[Condition B:]  $\quad \displaystyle p\to \infty,~\gama\ge \frac{c}{p} ~ \hbox{ and } ~c\ln \frac{1}{c}+(1-c)\ln \frac{1}{1-c}\ge \gama+ \gama\ln\frac{c+\gama}{\gama}+\delO$.
\end{description}
{\bf Remark~1}    {\bf Condition  A}  and  {\bf Condition  B} cover a wide range of $c$ and $\gama$, as illustrated below.
\begin{itemize}
\item[({\bf A})] For  each given prime $p$ and $c\le  \frac{p-1}{p+1}$, we have $1-c\ge \frac{1+c}{p}$ and hence 
\[
(1-c)\ln \frac{1}{1-c}> \frac{1+c}{p}\ln(1+c).
\]
Thus, {\bf Condition  A}  holds for all $c,\gama$ which satisfy
\begin{align*}
0<c\le  \frac{p-1}{p+1}\quad  \hbox{ and }~0\le \gama\le \frac{p-1}{p\ln(2p)}c\ln\frac{1}{c}.
\end{align*}
\item[({\bf B})] Suppose $q=p^a$ for some constant $a$. For example, $a=1$ corresponds to the {\em prime field}. Then $q\to \infty \Leftrightarrow p\to \infty$.  It is easy to see that {\bf Condition  B}  holds for any constant  $c\in (0,1)$ and all sufficiently small positive $\gama$. 
\end{itemize}
Define
\begin{align}
\mu_m(r)=\sum_{j=0}^{m}(-1)^j{\qb-r\choose j}q^{-j}. \label{eq:mu}
\end{align}
\begin{thm}\label{thm:Asy}  Let  $Q\in \Mcal_t$,  $\displaystyle D:=\{x\in \fq:Q(x)\ne 0\}$, $n:=|D|$ and  $\ell\ge 1$.  Suppose  either {\bf Condition  A} or {\bf Condition  B} holds. 
Then,  as  $k\to \infty$,  we have,  uniformly for   $0\le r\le k+t+\ell$ and $\eps\in \Ecal^{\ell,Q}$,
\begin{align}
\prob(Y_k(\eps)=r)&\sim \mu_{k+t+\ell-r}(r){\qb\choose r}q^{-r}. \label{eq:Asy}
\end{align}
\end{thm}

Recall that the (standard) Reed-Solomon code $\RS_{q,k}$ consists of the codewords $(g(x):x\in \fq)$ where  $g$ is a polynomial over $\fq$ of degree less than $k$.  When $Q=1$, we recall \cite{ChenMu07,LiWan20} that $q-Y_k(\langle f\rangle)$ is the distance between a received word $f\in \Mcal_{k+\ell}$ and a random codeword in $\RS_{q,k}$. Let $N(f,r)$ be the number of codewords in $\RS_{q,k}$ which are at distance $q-r$  from a received word $f$.  Thus, $ N(f,r)=q^k\prob(Y_k(\langle f\rangle)=r)$.  Setting $t=0$ in Theorem~\ref{thm:Asy}, we immediately obtain the following result.
\begin{cor}\label{cor:RS}  Suppose  either {\bf Condition  A} or {\bf Condition  B} holds.  Then, as  $k\to \infty$,  we have,  uniformly for   $0\le r\le k+\ell$ and $f\in \Mcal_{k+\ell}$,  
\[
N(f,r)\sim \mu_{k+\ell-r}(r){q\choose r}q^{k-r}.
\]
\end{cor}

\noindent {\bf Remark~2}  For  the prime field (i.e. $q=p$), Li and Wan \cite[Corollary~1.9]{LiWan20}  derived an asymptotic expression of $N(f,r)$ when the parameters satisfy the conditions: $k=cp$, $\ell=p^{\delta}$, $r=k+p^{\lambda}$, where $c\in (0,1)$, $\delta\in (0,1/4)$ and $\lambda\in (0,\delta)$ are all independent of $p$.  As commented in Remark~1, {\bf Condition B} is satisfied by any constant $c\in (0,1)$ and all sufficiently small positive $\gama$.  Thus,  Corollary~\ref{cor:RS}  covers $\ell$ up to $\sqrt{q}$ and all $r$, which significantly extends  the range of $\ell$ and $r$ covered by \cite[Corollary~1.9]{LiWan20}. 

Recall  that a received word represented by $f\in \Mcal_{k+\ell}$ is called a {\em deep-hole} if  $N(f,r)=0$ for all $r\ge k+1$, and is {\em  ordinary} if $N(f,k+\ell)>0$ (see, e.g.,  \cite{ChenMu07, GaoLi22,XuHong21}).  The well-known {\em Deep-Hole Conjecture} by Cheng and Murray  \cite{ChenMu07} states that  there is no  deep-hole when $\ell\ge 1$. Corollary~\ref{cor:RS} immediately implies the following result, which further supports the {\em Deep-Hole Conjecture}.
\begin{cor} 
Suppose  either {\bf Condition  A} or {\bf Condition  B} holds.  Then for sufficiently large $q$, every $f\in \Mcal_{k+\ell}$ is ordinary. Thus,  the  {\em Deep-Hole Conjecture} holds under either {\bf Condition  A} or {\bf Condition  B}. 
\end{cor}

\medskip

The remaining paper is organized as follows.
In Section~2 we recall some preliminary results about Hayes equivalence, Weil bounds, and sieve formulas. In Section~3 we extend the generating function approach from \cite{Gao21,GaoLi22} to general Hayes equivalence classes and give a proof of Theorem~1. Section~4 provides detailed error estimates which are needed for the proof of Theorem~2.  Section~5 provides the proof of Theorem~2.  Section~6 concludes the paper.

\section{Preliminaries}

In this section we recall some basic results needed to prove our main theorems. Hayes' theory of equivalence was first introduced in  \cite{Hay65}.   For $Q\in \Mcal_t$, define
\begin{align*}
\Phi_j(Q)&=|\{g\in \Mcal_j: \gcd(g,Q)=1\}|.
\end{align*}
In the rest of the paper, we shall use Iverson\rq{}s bracket $\llbracket P  \rrbracket$ which has value 1 if the predicate $P$ is true and 0 otherwise. 

It is easy to see \cite{Coh05,EH91}  that
\begin{align}
\left| \Mcal_{k+t+\ell}(\eps)\right|&=q^k, \label{eq:epsSize} \\
\left|\Ecal^{\ell,Q}\right|&=q^{\ell}\Phi_t(Q). \label{eq:Eorder}
\end{align}

 Let  $\{P_i:i\in I\}$ be the set of distinct irreducible factors of $Q$, where $P_i\in \Mcal_{d_i}$. Then the classical sieve formula gives
\begin{align*}
\Phi_j(Q)=\sum_{S\subseteq I}(-1)^{|S|} \left\llbracket  \sum_{i\in S}d_i\le j \right\rrbracket q^{j-\sum_{i\in S}d_i},
\end{align*}
and consequently
\begin{align}
 q^j\left(1-\sum_{i\in I}q^{-d_i}\right) \le \Phi_j(Q)\le q^j.  \label{PhiBound}
\end{align}
For example, if $I=\{1\}$ then
\begin{align*}
\Phi_j(Q)&= q^j -\llbracket j\ge d_1 \rrbracket q^{j-d_1}.
\end{align*}
If  $I=\emptyset$, that is, $Q=1$, then 
\begin{align*}
\Phi_j(Q)&= q^j.
\end{align*}

Noting $|I|\le \sum_{i\in I}d_i \le t$ and using  \eqref{PhiBound},  we obtain
\begin{align}\label{PhiGeneral}
\Phi_j(Q)&= q^j\left(1+O(t/q)\right).
\end{align}

For typographical convenience, we shall omit the superscripts and simply use $\Ecal$ to denote the group $\Ecal^{\ell,Q}$ when there is no danger of confusion.

Let ${\hat \Ecal}$ denote the group of characters over $\Ecal$, and $\chi\in {\hat \Ecal}$ be a  nontrivial character. Define $\chi(f)=\chi(\langle f\rangle)$ if $f\in \Mcal$ and $\gcd(f,Q)= 1$. Also set $\chi(f)=0$ if $f\in \Mcal$ and $\gcd(f,Q)\ne 1$ (This is the so called Dirichlet character).
By \cite[Ex. 5.2 \#2]{EH91} (see also \cite[Theorem~1.3]{Hsu96} and the paragraph before \cite[eq.~(4)]{Hsu96}), for each nontrivial character $\chi\in {\hat \Ecal}$, the associated $L$-function
\begin{align*}
P(z,\chi):&=\sum_{f\in \Mcal}\chi( f)z^{\deg(f)}
\end{align*}
is a polynomial of degree at most $\ell+t-1$. The roots of $P(z,\chi)$ are either 1 or have modulus $1/\sqrt{q}$.
Moreover, there is at most one root which is equal to  1. Thus we can write
\begin{align*}
P(z,\chi)&=\prod_{j=1}^{\deg(P(z,\chi))}(1-z\rho_j),\\
\deg(P(z,\chi))&\le \ell+t-1,~|\rho_j|\le \sqrt{q}.
\end{align*}
It follows  that
\begin{align}\label{eq:Wbound}
\left|\sum_{f\in \Mcal_j}\chi(f)\right|&\le {t+\ell-1\choose j}q^{j/2}.
\end{align}

Let $\ex(Y)$ denote the expected value of a random variable $Y$.  The following well-known result expresses  the probabilities in terms of the factorial moments \cite[Corollary~11]{Bo85}.
\begin{prop}\label{prop:sieve} Let $Y$ be any random variable which takes values in $\{0,1,\ldots, M\}$. We have
\begin{align*}
\prob(Y=r)&=\sum_{j=r}^{M}(-1)^{j-r}{j\choose r}\ex\left({Y\choose j}\right).
\end{align*}
Moreover, for each $r\le m\le M$, we have
\begin{align}
\left|\prob(Y=r)- \sum_{j=r}^{m-1}(-1)^{j-r}{j\choose r}\ex\left({Y\choose j}\right)\right|
\le {m\choose r}\ex\left({Y\choose m}\right). \label{eq:Alt}
\end{align}
\end{prop}

The following inequality \cite[(5)]{Bo85} will also be used to estimate binomial numbers.
\begin{align}\label{eq:bin}
{M\choose m} &\ge
  \left(\frac{M}{2\pi m(M-m)}\right)^{1/2}\left(\frac{M}{m}\right)^{m}\left(\frac{M}{M-m}\right)^{M-m} e^{-1/6}.  &(0< m<M)
\end{align}

The following are some simple observations about $\mu_m(r)$ defined in \eqref{eq:mu}. Since
\[
\frac{{\qb-r\choose j+1}q^{-j-1}}{{\qb-r\choose j}q^{-j}}=\frac{n-r-j}{q(j+1)}\le \frac{1}{j+1},
\]
the terms in \eqref{eq:mu} are alternating and have strictly decreasing absolute values.  Thus, we have
\begin{align}
\mu_m(r)&\ge \mu_1(r)= \frac{q-n+r}{q}. \label{eq:Alt}
\end{align}
Using the binomial expansion of  $\left(1-\frac{1}{q}\right)^{n-r}$, we also obtain
\begin{align}
\left|\mu_m(r)-\left(1-\frac{1}{q}\right)^{n-r}\right|&\le {n-r\choose m+1}\left(\frac{1}{q}\right)^{m+1}\le \frac{1}{(m+1)!}\left(\frac{n-r}{q}\right)^{m+1}\le \frac{1}{(m+1)!}. \label{eq:mubin}
\end{align}

Finally we recall the ``coordinate-sieve" formula by Li and Wan \cite{LiWan10,LiWan20}. Let $\Scal_j$ be the symmetric group on $\{1,2,\ldots,j\}$. 
Let $X$ denote the set of all $j$-tuples of elements from the set $D$, $\Xb$ the set of all  $j$-tuples of distinct elements from $D$. For each $\tau\in \Scal_j$, define 
\begin{align*}
X_{\tau}&=\{(x_1,\ldots,x_j)\in X, x_i= x_k \hbox{ if $i$ and $k$ belong to the same cycle of $\tau$} \}.
\end{align*}
 Let $h$ be a complex-valued function defined on $X$, and define
\begin{align*}
\Hb=\sum_{(x_1,\ldots,x_j)\in \Xb}h(x_1,\ldots,x_j),~~H(\tau)=\sum_{(x_1,\ldots,x_j)\in {X_\tau}}h(x_1,\ldots,x_j).
\end{align*}
Theorem~3.1 in  \cite{LiWan20} states
\begin{align}\label{eq:sieve}
\Hb=\sum_{\tau\in \Scal_j}(-1)^{j-l(\tau)}H(\tau),
\end{align}
where $l(\tau)$ denotes the number of cycles of $\tau$.

\section{Generating functions and proof of Theorem~1}

In this section, we use the generating function method developed in \cite{GKW22} to study the distribution of $Y_k(\eps)$.   It is also convenient to define $\langle f\rangle =0$ when $\gcd(f,Q)\ne 1$.

Define generating functions
\begin{align*}
F(z)&=\sum_{f\in \Mcal}\langle f\rangle z^{\deg(f)},\\
G(z,u)&=\sum_{f\in \Mcal} \langle f\rangle z^{\deg(f)}u^{r(f)},
\end{align*}
where $r(f)$ is the number of distinct zeros of $f$ that are in $D$.

We first prove the following result which is an extension of \cite[Prop.~2]{Gao21}:
\begin{prop}\label{prop:Gzu}
We have
\begin{align}
F(z)&=\sum_{d=0}^{t+\ell-1}\sum_{f\in \Mcal_d}\langle f\rangle z^d+\frac{1}{1-qz}z^{t+\ell}\sum_{\eps\in \Ecal}\eps,\label{eq:F}\\
G(z,u)&=F(z)\prod_{\alp\in D}(\langle 1\rangle+(u-1)z\langle x-\alp\rangle).\label{eq:G}
\end{align}
\end{prop}
\proof The proof is similar to that of \cite[Prop. 2]{Gao21};  see also \cite[p.27]{Gao21}.  Using \eqref{eq:epsSize}, we obtain
\begin{align*}
F(z)&=\sum_{d=0}^{t+\ell-1}\sum_{f\in \Mcal_d}\langle f\rangle z^d+\sum_{k\ge 0} q^k z^{k+t+\ell}\sum_{\eps\in \Ecal}\eps\\
&=\sum_{d=0}^{t+\ell-1}\sum_{f\in \Mcal_d}\langle f\rangle z^d+\frac{1}{1-qz}z^{t+\ell}\sum_{\eps\in \Ecal}\eps,
\end{align*}
which gives \eqref{eq:F}. 
Let $\Mcal'$ denote the subset of $\Mcal$ consisting of all  polynomials with no zeros in $D$. Then we have
\begin{align*}
F(z)&=\prod_{\alp\in D}\frac{1}{\langle 1\rangle-z\langle x-\alp\rangle}\sum_{g\in \Mcal'}z^{\deg (g)}\langle g\rangle, \nonumber\\
G(z,u)&=\prod_{\alp\in D}\left(\langle 1\rangle+u\sum_{j\ge 1}z^{j}\langle x-\alp\rangle^j\right)
\sum_{g\in \Mcal'}z^{\deg (g)}\langle g\rangle \nonumber\\
&=\left(\prod_{\alp\in D}\left(\frac{u}{\langle 1\rangle-z\langle x-\alp\rangle}+(1-u)\langle 1\rangle\right)\right)\left(F(z)\prod_{\alp\in D}(\langle 1\rangle-z\langle x-\alp\rangle)\right),
\end{align*}
which gives \eqref{eq:G}. \qed

As in \cite{GaoLi22}, let $D_j$ denote the set of all $j$-subsets of $D$. For $k+1\le j\le k+t+\ell$, define
\begin{align}
W_{j}(\eps)&=\sum_{g\in \Mcal_{k+t+\ell-j}}\sum_{S\in D_{j}}\left\llbracket \langle g\rangle \prod_{\alp\in S}\langle x-\alp\rangle=\eps \right\rrbracket. \label{eq:Wj}
\end{align}

For $\eps\in \Ecal$, we use $[z^d\eps]G(z,u)$ to denote the coefficient of $z^d\eps$ in $G(z,u)$.  Our next result gives all the factorial moments expressed in terms of $W_j(\eps)$. This extends the corresponding results for $Q=1$ in \cite{GaoLi22}.
\begin{thm}\label{thm:thm2} For each $Q\in \Mcal_t$ and $\eps\in \Ecal^{\ell,Q}$, we have
\begin{align}
\Ex\left({Y_k(\eps)\choose j}\right)&=\llbracket  j\le k\rrbracket {n\choose j}q^{-j}
+\llbracket k+1\le j\le k+t+\ell\rrbracket q^{-k}W_j(\eps).\label{eq:Moments}
\end{align}
\end{thm}
\proof Using Proposition~\ref{prop:Gzu}, we obtain
\begin{align}
G(z,u)&=\frac{1}{1-qz}z^{t+\ell}(1+(u-1)z)^n\sum_{\eps\in \Ecal}\eps\nonumber\\
&~~~~+\left(\sum_{j=0}^{t+\ell-1}z^{j}\sum_{g\in \Mcal_j}\langle g\rangle\right)\prod_{\alp\in D}\left(\langle 1\rangle+(u-1)z\langle x-\alp\rangle\right),\nonumber\\
\left[z^{k+t+\ell}\eps\right]G(z,u)&=\left[z^k \right]\frac{1}{1-qz}\left( 1+(u-1)z\right)^n\nonumber\\
&~~~~+\sum_{j=k+1}^{k+t+\ell}(u-1)^{j}[\eps]\sum_{g\in \Mcal_{k+t+\ell-j}}\sum_{S\in D_{j}}\langle g\rangle \prod_{\alp\in S}\langle x-\alp\rangle\nonumber\\
&=\sum_{j=0}^{k}{n\choose j}q^{k-j}(u-1)^{j}+\sum_{j=k+1}^{k+t+\ell}W_{j}(\eps)(u-1)^{j}.\label{eq:Gu}
\end{align}

Using \eqref{eq:Gu}, we obtain the following probability generating function of $Y_k(\eps)$:
\begin{align*}
p_{k}(u,\eps)&=q^{-k}\left[z^{k+t+\ell}\eps\right]G(z,u)\\
&=\sum_{j=0}^{k} q^{-j}{n\choose j}(u-1)^j
+q^{-k}\sum_{j=k+1}^{k+t+\ell}W_j(\eps)(u-1)^{j}.
\end{align*}
Hence 
\begin{align*}
\Ex\left({Y_k(\eps)\choose j}\right)&=\frac{1}{j!}\frac{d^j}{(du)^j}p_{k}(u,\eps)\Bigr|_{u=1} \\
&=\llbracket j\le k\rrbracket {n\choose j}q^{-j}
+\llbracket k+1\le j\le k+t+\ell\rrbracket q^{-k}W_j(\eps),
\end{align*}
which is \eqref{eq:Moments}.
 \qed

\medskip

\noindent {\bf Proof of Theorem~1} Substituting $m=k$ into \eqref{eq:Alt} and using  \eqref{eq:Moments},  
 we obtain
\begin{align*}
\left|\prob(Y_k(\eps)=r)-\sum_{j=r}^{k-1}(-1)^{j+r}{j\choose r}{n\choose j}q^{-j} \right|
&\le {k\choose r}{n\choose k}q^{-k}.
\end{align*}
Changing the summation index $j:=j+r$, and using  \eqref{eq:mu} and
\begin{align}\label{eq:Binoid}
 {j\choose r}{n\choose j}&= {n\choose r}{n-r\choose j-r},
\end{align}
 we obtain
\begin{align*}
\left|\prob(Y_k(\eps)=r)-{n\choose r}q^{-r} \mu_{k-1-r}(r) \right|
&\le {n\choose r}{n-r\choose k-r}q^{-k}\le  \frac{1}{(k-r)!} {n\choose r}q^{-r}.
\end{align*}
Using   \eqref{eq:mubin},  we obtain
\begin{align*}
\left|\prob(Y_k(\eps)=r)-{n\choose r}\left(\frac{1}{q}\right)^{r}\left(1-\frac{1}{q}\right)^{n-r} \right|
&\le \frac{1}{(k-r)!}{n\choose r}q^{-r}+  \frac{1}{(k-r)!} {n\choose r}q^{-r}\\
&\le  \frac{2}{(k-r)!}{n\choose r}\left(\frac{1}{q}\right)^{r}.
\end{align*}
Part (a) follows by using the assumption $k-r\to \infty$ and noting
\[
\left(1-\frac{1}{q}\right)^{n-r}\ge \left(1-\frac{1}{q}\right)^{q}\ge \frac{1}{4} \quad  (\hbox{ when  } q\ge 2).
\]
Part (b) follows immediately from part~(a) by noting 
\[
{n\choose r}\sim \frac{1}{r!}n^r,
\]
when  $r=o(\sqrt{n})$. \qed

\section{Estimates of $W_j(\eps)$}

In \cite{GaoLi22}, we derived estimate for $W_j(\eps)$ when $Q=1$ and $D=\fq$. In this section we carry out more detailed estimate for $W_j(\eps)$ and for general $Q$. Throughout this section, we  assume 
\begin{align}  \label{eq:Ddef}
D=\{\alpha\in \fq:Q(\alpha)\ne 0\}.
\end{align}
Then we have
\begin{align}  \label{eq:nBound}
q-t\le n \le q .
\end{align}

As in \cite{LiWan20}, the sum in \eqref{eq:Wj} can be estimated using the ``coordinate-sieve" formula \eqref{eq:sieve}.  Recall that $S_j$ denotes the set of all permutations of  $1,2,\ldots,j$. For $\tau\in S_j$,  $l(\tau)$ denotes the total number of cycles of $\tau$.  
We shall also use  $l'(\tau)$ to denote the number of cycles of $\tau$ which are not multiples of $p$. The standard generating function argument  \cite{FS09,Gao21} gives
\begin{align*}
\sum_{j\ge 0}\frac{1}{j!}\sum_{\tau\in S_j}a^{l(\tau)}b^{l'(\tau)}z^j&=\exp\left(a\sum_{i\ge 1,p\mid i}z^i/i+ab\sum_{i\ge 1,p\nmid i}z^i/i\right)\\
&=(1-z)^{-ab}(1-z^p)^{-(a-ab)/p}.
\end{align*}
As in \cite{Gao21,GaoLi22}, we define
\begin{align}
A_j(a,b)&= \frac{1}{j!}\sum_{\tau\in S_j}a^{l(\tau)} b^{l'(\tau)}\label{eq:cycle}\\
&=[z^j]\left((1-z)^{-ab}(1-z^p)^{-(a-ab)/p}\right)\label{eq:Agen}\\
&=\sum_{0\le i\le j/p}{ab+j-ip-1\choose j-ip}{(a-ab)/p+i-1 \choose i}.\nonumber
\end{align}

We have the following estimate of $W_j(\eps)$.
\begin{prop}\label{prop:Wj} Let $\eps\in \Ecal$, $k+1\le j\le k+t+\ell$. Suppose  $\ell\ge 1$ and $\gama:=(t+\ell-1)\sqrt{q}/n\le 1$. Then
\begin{align}
\left|W_j(\eps)-\frac{\Phi_{k+t+\ell-j}(Q)}{\Phi_t(Q)}{\qb\choose j}q^{-\ell}\right|&\le \frac{|\Ecal|-1}{|\Ecal|}{t+\ell-1\choose t+\ell+k-j}q^{(t+\ell+k-j)/2}A_j(\qb,\gama).\label{eq:Wjbound}
\end{align}
\end{prop}
\proof Let ${\hat \Ecal}$ denote the set of characters over the group $\Ecal$. Using \eqref{eq:Wj} and orthogonality of the characters, we obtain
\begin{align*}
W_{j}(\eps)&=\frac{1}{|\Ecal|}\sum_{g\in \Mcal_{k+t+\ell-j}}\sum_{S\in D_j}
\sum_{\chi\in {\hat \Ecal}}\chi\left( \eps^{-1}\langle g\rangle \prod_{\alpha\in S} \langle x-\alpha\rangle\right)\\
&=\frac{1}{|\Ecal|}\sum_{\chi\in {\hat \Ecal}}\chi(\eps^{-1})\left(\sum_{g\in \Mcal_{k+t+\ell-j}}\chi(g)\right)\sum_{S\in D_j}
\chi\left( \prod_{\alpha\in S} \langle x-\alpha\rangle\right).
\end{align*}

For the trivial character $\chi$, we have $\chi(\eps^{-1})=1$ and
\[
\sum_{g\in \Mcal_{k+t+\ell-j}}\chi(g)=\Phi_{k+t+\ell-j}(Q).
\]
It follows that
\begin{align}
W_{j}(\eps)&=\frac{\Phi_{k+t+\ell-j}(Q)}{|\Ecal|}{n\choose j}+\frac{1}{|\Ecal|}\sum_{\chi\ne 1}\chi(\eps^{-1})\left(\sum_{g\in \Mcal_{k+t+\ell-j}}\chi( g)\right) \sum_{S\in D_j}
\chi\left( \prod_{\alpha\in S} \langle x-\alpha\rangle\right).\label{eq:Wj2}
\end{align}
Applying \eqref{eq:Wbound}, we obtain
\begin{align}
\left|\sum_{g\in \Mcal_{k+t+\ell-j}}\chi( g)\right|&\le {t+\ell-1\choose t+\ell+k-j}q^{(k+t+\ell-j)/2}. \label{eq:eta}
\end{align}
Next, we use  \eqref{eq:sieve} to estimate  $\displaystyle \sum_{S\in D_j}
\chi\left( \prod_{\alpha\in S} \langle x-\alpha\rangle\right)$. 
For a permutation $\tau\in \Scal_j$, let $c_i$ be the number of cycles of length $i$ in $\tau$. Recall from the end of Section 2  that $X_{\tau}$ denotes the set of $j$-tuples of elements 
from $D$ which is constant in each cycle of $\tau$. As in \cite{LiWan20}, we set
\[
G_{\tau}:=\sum_{\vec{x}\in X_{\tau}}
\prod_{i=1}^{j} \chi(x-x_i).
\]
Then (Noting that $D$ contains all elements of $\fq$ which are not zeros of $Q$ and $\chi(x-\alpha)=0$ if $\alpha$ is a zero of $Q$.)
\begin{align*}
G_{\tau}&=\prod_{i} \left(\sum_{\alpha\in D}
\chi^i(x-\alpha)\right)^{c_i}\\
&=\left(\prod_{i,p\nmid i} \left(\sum_{\alpha\in \fq}\chi^{i}(x-\alpha)\right)^{c_i}\right)
\left(\prod_{i,p\mid i} \left(\sum_{\alpha\in \fq}\chi^{i}(x-\alpha)\right)^{c_i}\right).
\end{align*}
Since $\Ecal^{\ell,1}$ is a $p$-group, $\chi^i$ is nontrivial when $p\nmid i$.  It follows from  \eqref{eq:Wbound}  that 
\begin{align*}
\left|G_{\tau}\right|&\le
 \prod_{i,p\nmid i}(\gamma n)^{c_i}
\prod_{i,p\mid i} \qb^{c_i}
=(\gama n)^{l'(\tau)}\qb^{l(\tau)-l'(\tau)}=\qb^{l(\tau)}\gama^{l'(\tau)}.
\end{align*}
It follows from \eqref{eq:sieve} and \eqref{eq:cycle} that
\begin{align}\label{eq:linearbound}
\frac{1}{j!}\left|\sum_{\vec{x}\in \overline{X}}
\prod_{i=1}^j \chi(x-x_i)\right|
&\le  \frac{1}{j!}\sum_{\tau\in S_j}|G_{\tau}|\le  \frac{1}{j!}\sum_{\tau\in S_j}\qb^{l(\tau)}\gama^{l'(\tau)}=A_j(\qb,\gama).
\end{align}
Substituting \eqref{eq:eta} and \eqref{eq:linearbound} into \eqref{eq:Wj2}, and using \eqref{eq:Eorder} and 
\begin{align*}
\sum_{S\in D_{j}}\prod_{\alp\in S}\chi (x-\alp)
=\frac{1}{j!}\sum_{\vec{x}\in \overline{X}}
\prod_{i=1}^j \chi(x-x_i),
\end{align*}
we complete the proof. \qed

\begin{thm}\label{thm:thm4}   Let $\eps\in \Ecal^{\ell,Q}$ and assume $\ell\ge 1$.
We have
\begin{align}
&~~~\left|\prob(Y_k(\eps)=r)-\mu_{k-r}(r){\qb\choose r}q^{-r}-{n\choose r}q^{-(k+\ell)}\sum_{j=k+1}^{k+t+\ell}(-1)^{j-r}{n-r\choose j-r}\frac{\Phi_{k+t+\ell-j}(Q)}{\Phi_t(Q)}\right|\nonumber\\
&<q^{-k} \sum_{j=k+1}^{k+t+\ell} {j\choose r}{t+\ell-1\choose k+t+\ell-j}q^{(k+t+\ell-j)/2}A_j(\qb,\gama).\label{eq:thm4a}
\end{align}
\end{thm}
\proof   Using Theorem~4, inequalities (12) and  (26),  we obtain
\begin{align*}
&~~\left|\prob(Y_k(\eps)=r)-\sum_{j=r}^{k}(-1)^{j-r}{j\choose r}{\qb\choose j}q^{-j}-q^{-(k+\ell)}\sum_{j=k+1}^{k+t+\ell}(-1)^{j-r}{j\choose r}{\qb\choose j}\frac{\Phi_{k+t+\ell-j}(Q)}{\Phi_t(Q)}\right|\\
&\le \frac{|\Ecal|-1}{|\Ecal|}q^{-k}\sum_{j=k+1}^{k+t+\ell}\llbracket j\le n\rrbracket {j\choose r}{t+\ell-1\choose k+t+\ell-j}q^{(k+t+\ell-j)/2}A_j(\qb,\gama).
\end{align*}
Now \eqref{eq:thm4a} follows by using \eqref{eq:Binoid} and noting (as in the proof of Theorem~1)
\[
\sum_{j=r}^{k+t+\ell}(-1)^{j-r}{j\choose r}{\qb\choose j}q^{-j}=\mu_{k+t+\ell-r}(r){\qb\choose r}q^{-r}.
\]
 \qed

The following lemma provides the upper bounds for $A_j(n,\gama)$ which will be used in the proof of Theorem~\ref{thm:Asy}. Part~($a$)  is \cite[Lemma~1]{GaoLi22}.
\begin{lemma}\label{lem:Abound}
Let  $n\ge 1$. Then
\begin{itemize}
\item[(a)] For all $1\le j\le n$ and $b\in [0,1]$, we have
\begin{align}
\ln A_j(n,\gama)&\le \frac{j}{p}\ln \frac{n+j}{j}+\frac{n(1-\gama)}{p}\ln\frac{n+j}{n}
+n\gama\ln (2p). \label{eq:speciala}
\end{align}
\item[(b)] For all $1\le j\le 2pn\gama$ and $\gama\in (0,1]$,  we have
\begin{align}
\ln A_j(n,\gama) &\le j\ln \frac{n\gama+j}{j}+n\gama\ln \frac{n\gama+j}{n\gama}
+\frac{n(1-\gama)}{p}\ln 3.  \label{eq:specialb}
\end{align}
\end{itemize}
\end{lemma}
\proof From \eqref{eq:Agen}, we have
\begin{align*}
\sum_{ j\ge 0}A_j(n,\gama) z^j&=(1-z)^{-n\gama}(1-z^p)^{-n(1-\gama)/p}.
\end{align*}
It follows that, for $0<y<1$,
\begin{align}
A_j(n,\gama)&\le y^{-j}(1-y)^{-n\gama}(1-y^p)^{-n(1-\gama)/p}, \nonumber\\
\ln A_j(n,\gama) &\le j\ln y^{-1}+n\gama\ln(1-y)^{-1}+\frac{n(1-\gama)}{p}\ln\left(1-y^p\right)^{-1}. \label{eq:Abound1}
\end{align}
To minimize the above upper bound, we choose $y$  near the solution to the following saddle point equation
\begin{align}\label{eq:saddle}
-\frac{j}{y}+n\gama\frac{1}{1-y}+n(1-\gama)\frac{y^{p-1}}{1-y^p}=0,~~
i.e.,~~~\gama\frac{y}{1-y}+(1-\gama)\frac{y^p}{1-y^p}=\frac{j}{n}.
\end{align}

\noindent Part~(a) was proved in \cite{GaoLi22}. Since a similar argument will also be used for part~(b), we repeat it here. When $\gama$ is near $0$, we may  drop the term $\gama \frac{y}{1-y}$ in \eqref{eq:saddle}  to obtain the following approximate solution :
\begin{align*}
y=\left(\frac{j}{n+j}\right)^{1/p}.
 \end{align*}
Substituting this into \eqref{eq:Abound1}, we obtain
\begin{align}
\ln A_j(n,\gama)&\le \frac{j}{p}\ln \frac{n+j}{j}+\frac{n(1-\gama)}{p}\ln\left(\frac{n+j}{n}\right)
+n\gama\ln\left(1- \left(\frac{j}{n+j}\right)^{1/p}\right)^{-1}.\label{eq:general}
\end{align}
Since $(\ln 2)/p\le (\ln 2)/2<0.5$, we have
\begin{align*}
 \left(\frac{1}{2}\right)^{1/p}= \exp\left(\frac{-\ln 2}{p}\right)\le 1-\left(\frac{3}{4}\right)\left(\frac{\ln 2}{p}\right),
\end{align*}
and consequently
\begin{align}
\left(\frac{j}{n+j}\right)^{1/p}&\le \left(\frac{1}{2}\right)^{1/p}\le 1-\frac{3\ln 2}{4p},\nonumber\\
\ln\left(1- \left(\frac{j}{n+j}\right)^{1/p}\right)^{-1}&\le \ln \left(\frac{3\ln 2}{4p}\right)^{-1}
\le \ln (2p). \label{eq:bound}
\end{align}
Substituting \eqref{eq:bound} into \eqref{eq:general}, we obtain \eqref{eq:speciala}.

\noindent (b)   When  $p$ is large, $y^p$ is near zero and we drop the term $\displaystyle  (1-\gama)\frac{y^p}{1-y^p}$ in  \eqref{eq:saddle} to obtain the following approximate solution:
\[
y=\frac{j}{n\gama+j}.
\]
Substituting this into \eqref{eq:Abound1}, we obtain
\begin{align}
\ln A_j(n,\gama)&\le j\ln \frac{n\gama+j}{j}+n\gama\ln\left(\frac{n\gama+j}{n\gama}\right)
+\frac{n(1-\gama)}{p}\ln\left(1- \left(\frac{j}{n\gama+j}\right)^{p}\right)^{-1}.\label{eq:general2}
\end{align}
Using the assumption $\displaystyle  j\le 2pn\gama$, we obtain
\begin{align*}
 \left(\frac{j}{n\gama+j}\right)^{p}&\le  \left(1+\frac{1}{2p}\right)^{-p}\le \frac{2}{3}, \\
\ln\left(1- \left(\frac{j}{n\gama+j}\right)^{p}\right)^{-1}&\le  \ln 3.
\end{align*}
It follows from \eqref{eq:general2} that
\begin{align*}
\ln A_j(n,\gama)&\le j\ln \frac{n\gama+j}{j}+n\gama\ln\left(\frac{n\gama+j}{n\gama}\right)
+\frac{n(1-\gama)}{p}\ln 3.
\end{align*}
which is \eqref{eq:specialb}. \qed

\section{Proof of Theorem~2 }

{\bf Proof of Theorem~2:}   
Since the case $k-r\to \infty$ has already been covered by Theorem~1, we assume $r\ge k-\ln k$ below. Thus, we have  
\begin{align}\label{eq:rkBound}
t+\ell=O(\sqrt{n}),\quad q=n+O(\sqrt{n}), \quad \hbox{and }~r-k&=O(\sqrt{n}).
\end{align}
 It follows from \eqref{eq:Alt} that
$
 \mu_m(r)\ge c+O(1/\sqrt{n}),
$
which is bounded away from 0.
It follows from  \eqref{PhiGeneral} that
\begin{align}
&~~\left| q^{-k-\ell}\sum_{j=k+1}^{k+t+\ell}(-1)^{j-r}{n-r\choose j-r}\frac{\Phi_{k+t+\ell-j}(Q)}{\Phi_{t}(Q)}-\sum_{j=k+1}^{k+t+\ell}(-1)^{j-r}{n-r\choose j-r}q^{-j}\right|\nonumber \\
&= O\left(q^{-1/2}\right)\sum_{j=k+1}^{k+t+\ell}{n-r\choose j-r}q^{-j}\nonumber \\
&= o\left(q^{-r}\left(1+q^{-1}\right)^q\right) \nonumber \\
&= o\left(q^{-r}\right). \nonumber
\end{align}
By Theorem~\ref{thm:thm4},  asymptotic formula \eqref{eq:Asy} holds if
\begin{align*}
\sum_{j=k+1}^{k+t+\ell}{j\choose r}{t+\ell-1\choose k+t+\ell-j}q^{(t+\ell-k-j)/2}A_j(n,\gamma)=o\left({n\choose r}q^{-r}\right).
\end{align*}
Since  $t+\ell=O(\sqrt{q})$,  the above condition is implied by
\begin{align}\label{eq:Pbound}
\ln A_{j}(n,\gamma)+\ln {j\choose r}+\ln {t+\ell-1\choose t+\ell+k-r}-\ln {n\choose r}+\frac{t+\ell-k-j+1+2r}{2}\ln q\to -\infty.
\end{align}
Set
\begin{align*}
\delta&:=\frac{j}{\qb}-c=\frac{j}{\qb}-\frac{k}{q}, \quad (r\le j\le k+t+\ell),\\
\delta'&:=\frac{r}{\qb}-c=\frac{r}{\qb}-\frac{k}{q}.
\end{align*}
By \eqref{eq:rkBound},  we see that both $\delta$ and $\delta'$ are of the order $O\left(1/\sqrt{n}\right)$.
Using \eqref{eq:bin} and
\[
\ln(2\pi x(1-x))\le \ln (\pi/2)<1,~~(0<x<1)
\]
we obtain
\begin{align}
\ln {\qb\choose r}
&\ge \qb\left((c+\delta')\ln \frac{1}{c+\delta'}+(1-c-\delta')\ln \frac{1}{1-c-\delta'}\right)-\frac{1}{2}\ln \qb-\frac{2}{3}\nonumber\\
%& \ge \qb\left(c\ln \frac{1}{c}-\delta'+(1-c-\delta')\ln \frac{1}{1-c}\right)-\frac{1}{2}\ln \qb-\frac{2}{3}\nonumber\\
& \ge \qb\left(c\ln \frac{1}{c}+(1-c)\ln \frac{1}{1-c}\right)+O(\sqrt{n}). \label{eq:bin1}
\end{align}
Noting
\begin{align*}
\ln {j\choose r}&=\ln {j\choose j-r}\le (j-r)\ln j =O( \sqrt{n}\ln n),\\
\ln {t+\ell-1\choose t+\ell+k-r}&\le \ln 2^{t+\ell-1} = O( \sqrt{n}),
\end{align*}
we obtain
\begin{align}
&~~\ln {j\choose r}+\ln {t+\ell-1\choose t+\ell+k-r}-\ln {n\choose r}+\frac{t+\ell-k-j+1+2r}{2}\ln q \nonumber\\
&\le   -\qb\left(c\ln \frac{1}{c}+(1-c)\ln\frac{1}{1-c}\right)+O(\sqrt{n}\ln n). \label{Bound41}
\end{align}
We first consider {\bf Condition A}. Using \eqref{eq:speciala}, we obtain
\begin{align}
\ln A_{j}(n,\gamma)&\le \qb\left(\frac{c+\delta}{p}\ln \frac{1+c+\delta}{c+\delta}+\frac{1}{p}\ln\left(1+c+\delta\right)
+\gama\ln(2p)\right)\nonumber\\
&\le \qb\left(\frac{c}{p}\ln \frac{1+c}{c}+\frac{1}{p}\ln(1+c)
+\gama\ln(2p)\right)+O\left(\sqrt{n}\right).\nonumber
\end{align}
Thus
\begin{align*}
&~~~\ln A_{j}(n,\gamma)+\ln {j\choose r}+\ln {t+\ell-1\choose t+\ell+k-r}-\ln {n\choose r}+\frac{t+\ell-k-j+1+2r}{2}\ln q\\
&\le   -\qb\left(c\ln \frac{1}{c}+(1-c)\ln \frac{1}{1-c}-\frac{c}{p}\ln \frac{1+c}{c}-\frac{1}{p}\ln(1+c)-\gama\ln(2p)
+O\left(\frac{\ln n}{\sqrt{n}} \right) \right),
\end{align*}
which implies \eqref{eq:Pbound} under  {\bf Condition A}.

\medskip

Now we consider {\bf Condition B}. Using \eqref{eq:specialb} and 
\[
c\ln \left(1+\frac{\gama}{c}\right)\le \gamma,
\]
 we obtain
\begin{align*}
\ln A_j(n,\gamma)&\le j\ln \frac{\gamma n+j}{j}+\gama n\ln \frac{\gama n+j}{\gama n}+\frac{n(1-\gama)}{p} \ln 3\nonumber\\
&\le n\left((c+\delta)\ln \frac{c+\delta+\gama}{c+\delta}+\gama\ln \frac{c+\delta+\gama}{\gama}+\frac{1-\gama}{p} \ln 3\right)\\
&\le n\left(\gama+\gama\ln \frac{c+\gama}{\gama}+\frac{1}{p} \ln 3+O\left(\frac{1}{\sqrt{n}}\right)\right).
\end{align*}
It follows from  \eqref{Bound41} that
\begin{align*}
&~~~\ln A_{j}(n,\gamma)+\ln {j\choose r}+\ln {t+\ell-1\choose t+\ell+k-r}-\ln {n\choose r}+\frac{t+\ell-k-j+1+2r}{2}\ln q\\
&\le  - \qb\left(c\ln \frac{1}{c}+(1-c)\ln \frac{1}{1-c}-\gama-\gama\ln \frac{c+\gam}{\gama}
+O\left(\frac{\ln n}{\sqrt{n}}+\frac{1}{p}\right) \right),
\end{align*}
which  implies \eqref{eq:Pbound}  under {\bf Condition B}.   This completes the proof of Theorem~\ref{thm:Asy}. \qed

\medskip

\section{Conclusion}

For any given  $D\subseteq \fq$, we proved that the number of zeros in $D$  of a random monic polynomial over $\fq$ in a given Hayes\rq{}s equivalence class is asymptotically Poisson with mean $|D|/q$. We used the generating functions defined on the group algebra of equivalence classes, Weil's bounds on the corresponding character sums, and Li-Wan's coordinate-sieve formula.  It will be of considerable interest if the problem can be expressed as the sum of nearly independent random variables so that the Chen-Stein method can be applied to obtain the same results.  Asymptotic formulas are derived when the relevant parameters satisfy some simple inequalities.   Our asymptotic formulas extend  those in \cite{LiWan20,GaoLi22} about  the distance distribution in  Reed-Solomon codes, and  imply that all words for a large family of Reed-Solomon codes are ordinary, which further supports the {\em deep-hole} conjecture.

\end{document}